\documentclass[final,leqno,onefignum]{siamltex704}

\usepackage{amsfonts,latexsym,amsmath,enumerate,xcolor}
\usepackage[normalem]{ulem}
\usepackage{pstricks,epsfig}

\usepackage{tikz}
\definecolor{ble}{RGB}{66 66 111}
\definecolor{turk}{RGB}{112 147 219}
\definecolor{rosso}{RGB}{205 38 38}
\newrgbcolor{ross}{1 0 0}
\newrgbcolor{bleu}{66 66 111}

\newtheorem{remark}[theorem]{Remark}
\newtheorem{assumption}[theorem]{Assumption}

\newcommand{\dd}{\mathrm{d}}
\newcommand{\ii}{\mathrm{i}}

\newcommand{\ignore}[1]{}

\newcommand{\notmid}{\mid\kern-0.5em\not\kern0.5em}

\newcommand{\mstrut}[1]{\mbox{\rule{0mm}{#1}}}

\newcommand{\ba}{\begin{array}}
\newcommand{\ea}{\end{array}}

\title{Lipschitz stability of an inverse boundary value problem for a
  Schr\"{o}dinger type equation}

\author{Elena Beretta\thanks{Dipartimento di Matematica "Guido Castelnuovo"
Universita' di Roma "La Sapienza", Roma, Italy ({\tt
      beretta@mat.uniroma1.it})}
\and
Maarten V. de Hoop\thanks{Center for Computational and Applied
  Mathemematics, Purdue University, West Lafayette, IN 47907 ({\tt
    mdehoop@purdue.edu}).}
\and
Lingyun Qiu\thanks{Center for Computational and Applied Mathemematics,
  Purdue University, West Lafayette, IN 47907 ({\tt qiu@purdue.edu}).}
}

\begin{document}
\maketitle
\begin{abstract}
In this paper we study the inverse boundary value problem of determining the potential in the Schr\"{o}dinger equation from the knowledge of the Dirichlet-to-Neumann map, which is commonly accepted as an ill-posed problem in the sense that, under general settings, the optimal stability estimate is of logarithmic type. In this work, a Lipschitz type stability is established assuming a priori that the potential is piecewise constant with a bounded known number of unknown values.
\end{abstract}
\section{Introduction}
\label{sec:1}

In this paper, we investigate the stability for the inverse boundary
value problem of a Schr\"{o}dinger equation with complex potential,
$q(x)$ say. This encompasses the Helmholtz equation with attenuation, when $q(x) =
\omega^2 c^{-2}(x)$, where $c$ denotes the speed of propagation and
$\omega$ is the frequency, which can be complex. In fact, the imaginary part of $\omega c^{-1}(x)$ characterizes the attenuation of waves in the medium.

We begin with formulating the direct problem. Let $u \in H^1(\Omega)$
be the weak solution to the boundary value problem,
\begin{equation}\label{Helmholtz}
\left\{
\begin{array}{rl}
   (-\Delta + q(x)) u = & 0 ,\quad x \in \Omega ,
\\
                    u = & g ,\quad x\in \partial \Omega ,
\end{array}
\right.
\end{equation}
where $\Omega \subset \mathbb{R}^n$, $n \ge 2$ is a bounded connected
domain, $q \in L^\infty(\Omega)$ is a complex-valued function and $g$
is prescribed in the trace space $H^{1/2}(\Omega)$. The
Dirichlet-to-Neumann map is the operator $\Lambda_q :\ H^{1/2}(\Omega)
\rightarrow H^{-1/2}(\Omega)$ given by
\begin{equation}\label{DN_map}
   g \to \Lambda_q g
     = \left.\frac{\partial u}{\partial \nu}\right|_{\partial\Omega} ,
\end{equation}
where $\nu$ is the exterior unit normal vector to
$\partial\Omega$.

The inverse problem that we consider, consists in determining $q$ when
$\Lambda_q$ is known. This problem arises in geophysics, for example,
in reflection seismology assuming a description in terms of
time-harmonic scalar waves. The topic of this paper is the issue of
continuous dependence of $q$ from the Dirichlet-to-Neumann map
$\Lambda_q $. The continuous dependence is of fundamental
importance for the robustness of any reconstruction, as well as
for the development of convergent iterative reconstruction procedures
starting not too far from the solution (cf. \cite{Hoop2012}).
More precisely, it has been proved that Landweber iteration
reconstruction methods converge if the continuous dependence for the
inverse problem is of H\"{o}lder or Lipschitz type.

From the work of \cite{Mandache2001}, it is evident that for arbitrary
potentials $q$, Lipschitz stability cannot hold. Motivated by, and
following analogous results in electrical impedance tomography (EIT,
cf. \cite{Alessandrini2005,Beretta2011}), here we study
conditional stability when a-priori information on $q$ is assumed. We
consider models with discontinuous potentials to accommodate realistic
reflectors. Specifically, we consider the space spanned by linear
combinations of $N$ characteristic functions. More precisely we
consider potentials of the form
\[
   q(x) = \sum_{j = 1}^N q_j \chi_{D_j}(x) ,
\]
where $q_j, j=1,\dots N$ are unknown complex numbers and $D_j$ are
known open Lipschitz sets in $\mathbb{R}^n$.  Moreover, we consider
the case of partial boundary data, that is, we can restrict the
collection of measurements to only a part of the boundary. We refer to
\cite{Uhlmann2009} for a review of recent uniqueness results. Here, we
prove Lipschitz stability with a uniform constant, which depends on
$N$ and on the other a-priori parameters of the problem. We will show
that the Lipschitz constant grows exponentially with the dimension,
$N$, of the space of potentials. The method of proof follows the ideas
introduced in Alessandrini and Vessella and relies on quantitative
estimates of unique continuation of solutions to elliptic systems and
on the use of singular solutions and of their asymptotic behaviour
near the discontinuity interfaces.  Compared to the case of the real
or complex conductivity equation in the case of the Schr\"{o}dinger
equation we are able to derive our result relaxing the assumptions of
regularity on $\partial D_j$ that are assumed to be Lipschitz.
Furthermore, taking advantage of the regularity of solutions and of
its gradient inside the domain $\Omega$ we find a better dependence of
the stability constant on $N$.

The outline of the paper is as follows. In the next section we state
all the assumptions and the main result. In Section~\ref{sec:3}, we
give a summary of known regularity results connected to
Schr\"{o}dinger equation with complex potential, and some preparatory
lemmas concerning the existence and asymptotic behaviour of singular
solutions. Section~\ref{sec:4} contains the proof of our main
theorem. We first show the proof for $n=3$ and then modify it to the other cases. For the structure of the main proof we characterize the rate of blow-up of the singular functions finding lower and upper bounds in terms of the distance of the singularity from the interface of the subdomains. More precisely, to derive our main result we first establish that the singular function satisfies a lower bound in terms of the distance of the singularity from the interface. Secondly, by using quantitative estimates of propagation of smallness we derive also an upper bound for the singular function. Last but not least, we make use the value of a bounded non-decreasing function at some particular point to prove that either the result of the main theorem can be deduced directly or a recursive inequality (\ref{induc-equ2}) must hold true. The recursive inequality also leads to the desired result. In Section~\ref{sec:5} we demonstrate by an example that the
Lipschitz constant grows exponentially with the dimension of the space
of potentials. This example is constructed from its analogue in
electrical impedance tomography \cite{Rondi2006}.

\section{Main result}\label{sec:2}

\subsection{Notation and definitions}\label{sec:2-1}

We denote by $n$ the space dimension. For every $x \in \mathbb{R}^n$,
we set $x=(x', x_n)$ where $x'\in \mathbb{R}^{n-1}$ for $n\ge 2$. With
$B_R(x)$, $B'_R(x')$ and $Q_R(x)$ we denote the open ball in
$\mathbb{R}^n$ centered at $x$ of radius $R$, the ball in
$\mathbb{R}^{n-1}$ centered at $x'$ of radius $R$, and the cylinder
$B'_R(x') \times (x_n - R, x_n + R)$, respectively. For simplicity of
notation, $B_R(0)$, $B'_R(0)$ and $Q_R(0)$ are denoted by $B_R$,
$B'_R$ and $Q_R$.

\medskip\medskip

\begin{definition}\label{def:Lip-portion}
Let $\Omega$ be a bounded domain in $\mathbb{R}^n$. We say that a
portion $\Sigma$ of $\partial\Omega$ is of Lipschitz class with
constants $r_0, L>0$ if, for any $P\in \Sigma$, there exists a rigid
transformation of coordinates such that $P = 0$ and
\[
   \Omega \cap Q_{r_0}
         = \{(x',x_n)\in Q_{r_0} \mid x_n > \phi(x') \}
\]
where $\phi$ is a Lipschitz continuous function on $B_{r_0}'$ with
$\phi(0)=0$ and
\[
   \|\phi\|_{C^{0,1}(B_{r_0}')} \le L.
\]
We shall say that $\Omega$ is of Lipschitz class with constants $r_0$
and $L$, if $\partial \Omega$ is of Lipschitz class with the same
constants.
\end{definition}

\medskip\medskip

\begin{definition}
Let $\Omega$ be a bounded open subset of $\mathbb{R}^n$ and of
Lipschitz class and $\Sigma$ be a open portion of $\partial
\Omega$. We define $H^{1/2}_{co}(\Sigma)$ as
\[
   H^{1/2}_{co}(\Sigma)
   = \{g \in H^{1/2}(\partial \Omega)
                 \mid \supp g \subset \Sigma\}
\]
and $H^{-1/2}_{co}(\Sigma)$ as the topological dual of
$H^{1/2}_{co}(\Sigma)$; we denote by $\langle \cdot, \cdot \rangle$
the dual pairing between $H^{1/2}_{co}(\Sigma)$ and
$H^{-1/2}_{co}(\Sigma)$.

\end{definition}

\medskip\medskip

\begin{definition}\label{def:loc-DtoN}
Let $\Omega$ be a bounded open subset of $\mathbb{R}^n$ and of
Lipschitz class, $\Sigma$ be a open portion of $\partial \Omega$ and
$q\in L^{\infty}(\Omega)$. Assume that $0$ is not an eigenvalue of
$(-\Delta + q)$ with Dirichlet boundary conditions in $\Omega$, i.e.,
\[
   \{u \in H^1_0(\Omega) \mid (-\Delta + q)u = 0\} = \{0\} .
\]
For any $g \in H^{1/2}_{co}(\Sigma)$, let $u\in H^1(\Omega)$ be the
weak solution to the Dirichlet problem
\begin{equation}\label{Dirichlet-pro}
\left\{\begin{array}{rl}
   (-\Delta + q(x))u = & 0 ,\quad x\in \Omega ,
\\
   u = & g ,\quad x\in \partial \Omega .
\end{array}\right.
\end{equation}
We define the local Dirichlet-to-Neumann map $\Lambda_{q}^{(\Sigma)}$
as
\[
\begin{array}{rrl}
   \Lambda_q^{(\Sigma)}: & H^{1/2}_{co}(\Sigma)
          \rightarrow & H^{-1/2}_{co}(\Sigma)
\\
   & g \mapsto &
      \displaystyle{\left.\frac{\partial u}{
        \partial \nu}\right|_{\Sigma}} ,
\end{array}
\]
where $\nu$ is the exterior unit normal vector to $\partial\Omega$.
\end{definition}

\medskip\medskip

With $\Omega$ being a bounded open set, with $C^{0,1}$ boundary, the
set of the eigenvalues of $(-\Delta + q)$ with Dirichlet boundary
conditions is a discrete subset of $\mathbb{C}$, and hence can be
avoided.

We observe that $\Lambda_q^{(\Sigma)}$ can be identified with the
sesquilinear form on $H^{1/2}_{co}(\Sigma) \times
H^{1/2}_{co}(\Sigma)$, defined by
\[
   \langle \Lambda_q^{(\Sigma)}g, f \rangle =
   \int_{\Omega} (\nabla u \cdot \nabla \bar v
   + q u \bar v) \dd x, \quad \forall f,g \in H^{1/2}_{co}(\Sigma) ,
\]
where $u$ is the solution to $(\ref{Dirichlet-pro})$ and $v$ is any
function in $H^1(\Omega)$ such that $v \mid_{\partial \Omega} =
f$. This definition is independent of the choice of $v$: Let $v_1,
v_2$ be two different functions in $H^1(\Omega)$ such that $v_1
\mid_{\partial \Omega} = v_2 \mid_{\partial \Omega}= f$. Then, since
$w = v_1 - v_2\in H^1_0(\Omega)$, and $u$ is a solution, we have
\begin{multline*}
   \int_{\Omega} \left( \nabla u \cdot \nabla \bar v_1
      + q u \bar v_1 \right) \dd x
   - \int_{\Omega} \left( \nabla u \cdot \nabla \bar v_2
      + q u \bar v_2 \right) \dd x
\\
   = \int_{\Omega} \left( \nabla u \nabla \bar w
      + q u \bar w \right) \dd x = 0 ,
\hspace*{1.0cm}
\end{multline*}
using integration by parts. We denote by
$\|\cdot\|_{\mathcal{L}(H^{1/2}_{co}(\Sigma), H^{-1/2}_{co}(\Sigma))}$
the norm defined as
\[
   \|\Lambda_q^{(\Sigma)}\|_{\mathcal{L}(H^{1/2}_{co}(\Sigma),
           H^{-1/2}_{co}(\Sigma))} =
   \sup_{f,g \in H^{1/2}_{co}(\Sigma)} \{
   \langle \Lambda_q^{(\Sigma)}g, f \rangle \mid
   \|g\|_{H^{1/2}_{co}(\Sigma)} = \|f\|_{H^{1/2}_{co}(\Sigma)} = 1 \} .
\]

\subsection{Main assumptions}
\label{sec:2-2}

Our assumptions on $\Omega$ and $q(x)$ are

\medskip\medskip

\begin{assumption}\label{assumption_domain}
$\Omega \subset \mathbb{R}^n$ is a bounded domain satisfying
\[
   |\Omega| \le A
\]
Here and in the sequel $|\Omega|$ denotes the Lebesgue measure of
$\Omega$. We assume that $\partial \Omega$ is of Lipschitz class and
we fix an open portion $\Sigma$ of $\partial\Omega$ which is of
Lipschitz class with constants $r_0$ and $L$.
\end{assumption}

\medskip\medskip

\begin{assumption}\label{assumption_potential}
The complex-valued function $q(x)$ satisfies
\[
\| q \|_{L^{\infty}(\Omega)} \le B,
\]
where $B$ is a positive constant, and is of the form
\[
   q(x) = \sum_{j = 1}^N q_j \chi_{D_j}(x) ,
\]
where $q_j, j=1,\dots N$ are unknown complex numbers and $D_j$ are
known open sets in $\mathbb{R}^n$ which satisfy the following
assumption. Moreover, we assume that $0$ is not an eigenvalue of
$-(\Delta + q)$ with Dirichlet boundary conditions in $\Omega$.
\end{assumption}

\medskip\medskip

\begin{assumption}\label{assumption_piece}
The $D_j, j=1,\dots,N$, are connected and pairwise non-overlapping
open sets such that $\cup_{j=1}^N \overline{D}_j = \overline{\Omega}$
and $\partial D_j$ are of Lipschitz class. We also assume that there
exists one set, say $D_1$, such that $\partial D_1 \cap \partial
\Omega$ contains an open portion $\Sigma_1$ of Lipschitz class with
constants $r_0$ and $L$. For every $j \in \{2, \dots, N\}$ there exist
$j_1,\dots,j_M \in \{1,\dots, N\}$ such that
\[
   D_{j_1} = D_1, \quad D_{j_M} = D_j
\]
and, for every $k = 1, \dots, M$,
\[
\partial D_{j_{k - 1}} \cap \partial D_{j_k}
\]
contains a non-empty open portion $\Sigma_k$ of Lipschitz
class with constants $r_0$ and $L$ such that
\[
\begin{array}{c}
\Sigma_1 \subset \Sigma, \\
\Sigma_k \subset \Omega, \quad \forall k =2, \dots , M.
\end{array}
\]

Furthermore, there exists $P_k\in \Sigma_k$, at which $D_{k-1}$
satisfies the interior ball condition with radius $\frac{3r_0}{16}$,
and a rigid transformation of coordinates such that $P_k = 0$ and
\[
\begin{array}{rl}
\Sigma_k \cap Q_{r_0/3} = & \{x\in Q_{r_0 /3} \mid x_n = \phi_k(x')\},\\
D_{j_k} \cap Q_{r_0/3} = & \{x\in Q_{r_0 /3} \mid x_n > \phi_k(x')\},\\
D_{j_{k - 1}}  \cap Q_{r_0/3} = & \{x\in Q_{r_0 /3} \mid x_n < \phi_k(x')\},
\end{array}
\]
where $\phi_k$ is a $C^{0,1}$ function on $B_{r_0/3}'$ satisfying
\[
\phi_k(0) = 0
\]
and
\[
\|\phi_k\|_{C^{0,1}(B_{r_0/3}')} \le L.
\]
For simplicity, we call $D_{j_1}, \dots, D_{j_M} $ a chain of domains
connecting $D_1$ to $D_j$.
\end{assumption}

\vskip 12truemm

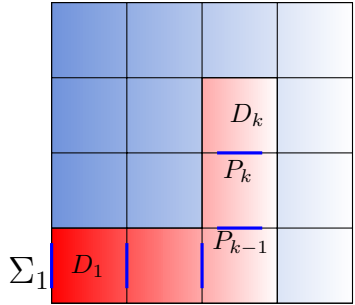
\begin{center}
\begin{tikzpicture}[scale=2]
\shadedraw[left color=turk] (0,0) rectangle (2,2);
\filldraw [left color=red] (0,0) -- (1.5,0) -- (1.5,1.5) --  (1,1.5) -- (1,0.5) -- (0,0.5) -- (0,0);
\draw (-0.35,0.2) node [right] {\Large $\Sigma_1$};
\draw (1.3,1.25) node { $D_{k}$};
\draw (0.24,0.25) node { $D_1$};
\draw (1.24,0.88) node {$P_k$};
\draw (1.26,0.4) node {$P_{k-1}$};
\draw[step=0.50] (0,0) grid (2,2);
\draw [very thick, color=blue] (0,0.1) -- (0,0.4);
\draw [very thick, color=blue] (0.5,0.1) -- (0.5,0.4);
\draw [very thick, color=blue] (1,0.1) -- (1,0.4);
\draw [very thick, color=blue] (1.1,0.5) -- (1.4,0.5);
\draw [very thick, color=blue] (1.1,1) -- (1.4,1);
\end{tikzpicture}
\end{center}
\medskip\medskip

In the further analysis, for simplicity of notation, we also use
the constant $r_1 = \frac{r_0}{16}$.

\subsection{Statement of the main result}\label{sec:2-3}

The main result of this paper is stated as follows.

\medskip\medskip

\begin{theorem}\label{thm-stability}
Let $\Omega$ satisfy Assumption~\ref{assumption_domain} and $q^{(k)},
k = 1, 2$ be two complex piecewise constant functions of the form
\[
   q^{(k)}(x) = \sum_{j = 1}^{N}
         q_j^{(k)} \chi_{D_j}(x) ,\quad k = 1, 2
\]
which satisfy Assumption~\ref{assumption_potential} and $D_j,
j=1,\dots, N$ satisfy Assumption~\ref{assumption_piece}. Then, there
exists a constant $C = C(n, r_0, L, A, B, N)$, such that
\begin{equation}\label{Lip_stab}
   \| q^{(1)} - q^{(2)}\|_{L^\infty(\Omega)} \le
       C \|\Lambda_1^{(\Sigma)}
   - \Lambda_2^{(\Sigma)}\|_{\mathcal{L}(H^{1/2}_{co}(\Sigma),
                H^{-1/2}_{co}(\Sigma))} ,
\end{equation}
where $\Lambda_k^{(\Sigma)} = \Lambda_{q^{(k)}}^{(\Sigma)}$ for
$k=1,2$.
\end{theorem}

\section{Preliminary results}
\label{sec:3}

In this section, we state some results which will be used in the proof of our main stability result.

\medskip\medskip

\begin{proposition}\label{p-energy}
Let $\Omega$ be a bounded Lipschitz domain in $\mathbb{R}^n$,  $q \in L^\infty(\Omega)$ complex valued potential , $f\in L^p(\Omega)$ and $g\in W^{2-\frac{1}{p},p}(\partial\Omega)$ with $1<p<\infty$. Assume that $0$ is not a Dirichlet eigenvalue for the operator $-\Delta+q$ in $\Omega$. Then there exists a unique solution $u\in W^{2,p}(\Omega)$ to the problem
\begin{equation}\label{Helmholtz-L_2}
\left\{
\begin{array}{rl}
(-\Delta + q(x))u = & f, \quad x\in \Omega, \\
u= & g, \quad x\in \partial \Omega,
\end{array}
\right.
\end{equation}
Moreover,
\begin{equation}\label{energy_est1}
\| u \|_{W^{2,p}(\Omega)} \le C \left( \| g \|_{W^{2-\frac{1}{p},p}(\partial \Omega)} + \| f \|_{L^{p}(\Omega)} \right)
\end{equation}
where $C$  depends on $n,\Omega$ and $\|q\|_{L^\infty(\Omega)}$.
\end{proposition}

\medskip\medskip

The proof is a consequence of the of existence of a $W^{2,p}(\Omega)$
function $w$ such that $w=g$ on $\partial\Omega$ and such that
$\|w\|_{W^{2,p}(\Omega)}\le C\|g\|_{W^{2-\frac{1}{p},p}(\partial\Omega)}$ and of the
Fredholm alternative; see for example Theorem 3.5.8 in Feldman and
Uhlmann's notes \cite{Feldman}). For reader's convenience, we also note the following Proposition~\ref{2-energy} without proof, which we use for the low dimension cases. 

\medskip\medskip

\begin{proposition}\label{2-energy}
Let $\Omega$ be a bounded Lipschitz domain in $\mathbb{R}^n$,  $q \in L^\infty(\Omega)$ complex valued potential , $f\in H^{-1}(\Omega)$ and $g\in H^{1/2}(\partial\Omega)$. Assume that $0$ is not a Dirichlet eigenvalue for the operator $-\Delta+q$ in $\Omega$. Then there exists a unique solution $u\in H^1(\Omega)$ to the equation (\ref{Helmholtz-L_2}).
Moreover,
\begin{equation}
\| u \|_{H^{1}(\Omega)} \le C \left( \| g \|_{H^{1/2}(\partial \Omega)} + \| f \|_{H^{-1}(\Omega)} \right)
\end{equation}
where $C$  depends on $n,\Omega$ and $\|q\|_{L^\infty(\Omega)}$.
\end{proposition}

\medskip\medskip

Our approach follows the one of Beretta and
Francini\cite{Beretta2011}, which is for the EIT problem with
complex conductivity, of constructing singular solutions and of
studying their asymptotic behavior when the singularity approaches the
interfaces $\Sigma_k$. This method was originally introduced by
Alessandrini and Vessella in the real-valued conductivity case
\cite{Alessandrini2005}. To construct singular solutions for the EIT
problems, the Green's function plays a crucial role. In our case, we
also use the Green's function to treat the case of high dimension ($n
\ge 4$) and a first order derivative of Green's function needs to be
used for lower dimension ($n = 2,3$). In the following propositions,
we discuss the existence and behavior of the Green's functions ($n \ge
4$) and a first order derivative of the Green's function ($n = 2,3$)
when $q$ satisfies Assumption~\ref{assumption_potential}. We are
especially interested in their asymptotic behavior near the $C^{0,1}$
interface $\Sigma_k$.

Before doing this, we need to extend our original domain. We consider
$\Sigma_1$ and recall that up to a rigid transformation of coordinates
we can assume that $P_1 = 0$ and
\[
   (\mathbb{R}^n \backslash \Omega) \cap B_{r_0}
       = \{(x', x_n) \in B_{r_0} \mid x_n < \phi(x')\}
\]
where $\phi$ is a Lipschitz function such that $\phi(0) = 0$ and $\|
\phi \|_{C^{0,1}(B_{r_0}')} \le L $. Then we extend $\Omega$ to
$\Omega_0=\Omega \cup D_0$ by adding an open set $D_0$ defined as
\[
   D_0 = \left\{ x \in (\mathbb{R}^n \backslash \Omega)
        \cap B_{r_0} \mid
     \left| x_n - \frac{r_0}{6} \right| < \frac{5}{6} r_0 ,
     \, |x_i|< \frac{2}{3}r_0, \, i=1,\dots, n-1 \right\} .
\]
It turns out that $\Omega_0$ is of Lipschitz class with constants
$\tfrac{r_0}{3}$ and $L_1$, where $L_1$ depends on $L$ only. We define
\[
   K_0 = \left\{ x \in D_0 \mid
           \operatorname{dist}(x,\Sigma_1) \ge \frac{r_0}{3} \right\}
\]
with $\operatorname{dist}(K_0, \partial \Omega) > \tfrac{r_0}{3}$. We
extend $q(x)$ defined on $\Omega$ by setting it equal to $1$ in
$D_0$. For simplicity of notation we still denote this extension by
$q(x)$.

We consider any subdomain in $\Omega$ and the chain of domains
connecting it to $D_1$. For simplicity let us rearrange the indices of
subdomains so that this chain corresponds to $D_0,D_1,\dots,D_M$, $M
\le N$. Let $S = \cup_{j = 0}^M \overline{D}_j$ and $K$ be a connected
subset of $S$ with Lipschitz boundary such that $\overline{K} \cap
\partial D_j = \Sigma_j \cup \Sigma_{j+1}$ for $j = 1,2, \dots, M$,
$K_0 \subset K$ and $\operatorname{dist}(K, \partial S \backslash \{\Sigma_{M+1} \cup
\Sigma_1\}) > \tfrac{r_0}{16}$.

\vskip 12truemm

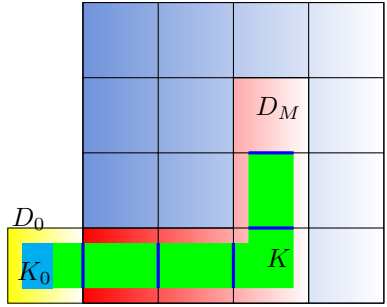
\begin{center}
\begin{tikzpicture}[scale=2]
\shadedraw[left color=turk] (0,0) rectangle (2,2);
\shadedraw [left color=yellow] (-0.5,0) rectangle (0,0.5);
\fill [cyan] (-0.4,0.1) rectangle (-0.2,0.4);
\draw (-0.5,0.2) node [right] {$K_0$};
\filldraw [left color=red] (0,0) -- (1.5,0) -- (1.5,1.5) --  (1,1.5) -- (1,0.5) -- (0,0.5) -- (0,0);
\fill [green] (-0.2,0.1) -- (1.4,0.1) -- (1.4,1) --  (1.1,1) -- (1.1,0.4) -- (-0.2,0.4);
\draw (1.3,1.3) node { $D_M$};
\draw (1.3,0.3) node { $K$};
\draw (-0.36,0.56) node { $D_0$};
\draw[step=0.50] (0,0) grid (2,2);
\draw [very thick, color=blue] (0,0.1) -- (0,0.4);
\draw [very thick, color=blue] (0.5,0.1) -- (0.5,0.4);
\draw [very thick, color=blue] (1,0.1) -- (1,0.4);
\draw [very thick, color=blue] (1.1,0.5) -- (1.4,0.5);
\draw [very thick, color=blue] (1.1,1) -- (1.4,1);
\end{tikzpicture}
\end{center}
\vskip 10truemm
In the following, we shall use $C$ to denote positive constants. The
value of the constants may change from line to line, but we shall
specify their dependence everywhere where they appear. For $n \ge 4$,
let $\Gamma$ denote the fundamental solution associated with the
Laplace operator. In the proof of Theorem~\ref{thm-stability}, we will need to estimate $G - \Gamma$ from above in terms of variable-interface distance $r$ to a power, which is smaller than the order of the singularity of $\Gamma$. Since, for high dimension cases($n \ge 6$), $\Gamma(\cdot,y)$ does not belong to $H^{-1}(\Omega)$, we need to employ $L^p$ estimate of the solutions here. Note that $\Gamma(\cdot,y)$ belongs to $L^p(\Omega)$ for any $1 \le p < \frac{n}{n-2}$. 

\medskip\medskip

\begin{proposition}\label{G_behavior}
Let the complex-valued function $q \in L^\infty(\Omega_0)$ satisfy
Assumption~\ref{assumption_potential} and $n \ge 4$. For $y \in
\Omega_0$, there exists a unique function $G(\cdot, y)$ continuous in
$\Omega_0 \backslash \{ y \}$ such that
\begin{equation}\label{fS}
   \int_{\Omega_0} \nabla G(\cdot,y) \nabla \phi
       + q G(\cdot,y) \phi = \phi(y) ,\quad
   \forall \phi \in C_0^\infty(\Omega) .
\end{equation}
Furthermore, we have that $G(x,y)$ is symmetric, that is,
\begin{equation}\label{symmetry}
   G(x,y) = G(y,x) ,\quad x,y \in \Omega_0 ,
\end{equation}
and the following estimates
\begin{equation}\label{singular_func_est}
   \begin{aligned}
   \| G(\cdot,y) \|_{L^2(\Omega_0 \backslash B_r(y))}
        \le & \,\, C |\ln r|^{\frac{1}{2}} ,\quad
   & r \le \frac{1}{2} \operatorname{dist}(y,\partial \Omega_0), \quad & n = 4
   \\
   \| G(\cdot,y) \|_{L^2(\Omega_0 \backslash B_r(y))}
        \le & \,\, C r^{2 -\frac{n}{2}} ,\quad
   & r \le \frac{1}{2} \operatorname{dist}(y,\partial \Omega_0), \quad & n \ge 5
   \end{aligned}
\end{equation}
and
\begin{equation}\label{est_green2}
\| G(\cdot, y) - \Gamma(\cdot,y) \|_{L^2(\Omega_0)} \le
\left\{
\begin{aligned}
  &\quad \quad C &&, \quad & 4 \le n\le 7 ,
  \\
  &|\ln(\operatorname{dist}(y, \cup_{j=1}^N\partial D_j))| &&, \quad & n=8 ,
  \\
  &{\operatorname{dist}(y, \cup_{j=1}^N\partial D_j)}^{4 - \frac{n}{2}} &&, \quad & n\ge 9,
\end{aligned}
\right.
\end{equation}
for $\operatorname{dist}(y,\partial \Omega_0)\ge \frac{r_0}{16}$ ,
hold true, where the constant $C$ depends on the constant in Proposition~\ref{p-energy}.
\end{proposition}

\medskip\medskip

\begin{proof}
Assume that $y$ belongs to some sub-domain $D_m$ which $q$ equals to a complex constant $q_m$ inside. Let $H(x,y)$ denote the outgoing fundamental solution of Helmholtz equation
 \[
(-\Delta + q_m)H(x,y) = \delta(x,y), \quad x\in \mathbb{R}^n,
 \]
 i.e.,
\begin{equation}\label{fundamental_Helmholtz}
 H(x,y) = \frac{ q_m^{(n-2)/4} H^{(1)}_{(n-2)/2}(q_m^{1/2}|x-y|)}{4 \ii (2\pi)^{(n-2)/2} |x-y|^{(n-2)/2}},
\end{equation}
 where $H^{(1)}_n$ denotes Hankel function of the first kind.
We consider $G(x,y) = H(x,y) +\omega (x,y)$, where $\omega$ solves
\begin{equation}
\left\{
\begin{array}{rl}
(-\Delta + q)\omega = & (q_m - q) H, \quad \mbox{ in }\Omega_0, \\
\omega = & - H, \quad \mbox{ on }\partial \Omega_0.
\end{array}
\right.
\end{equation}
Note that $q_m - q$ vanishes in $D_m$. Hence $(q_m-q)H$ belongs to $L^\infty(\Omega_0)$. By using the asymptotic behavior of the Hankel function near the origin \cite{Ruiz2002}, we obtain that
\[
|(q_m-q)H(x,y)| \le \left\{
\begin{aligned}
  &\quad \quad 0 &&, & \quad |x-y|\le \operatorname{dist}(y, \partial D_m),
  \\
  &C|x-y|^{2-n} &&, & \quad |x-y|> \operatorname{dist}(y, \partial D_m),
\end{aligned}
\right.
\]
for some positive constant $C$. We observe that the order of the singularity of $\omega(x,y)$ is always lower then the fundamental solution $H(x,y)$. To be more precise, by applying Proposition~\ref{p-energy} with $p = \frac{2n}{n+4}$ and Sobolev embedding theorem, we conclude that
\begin{multline}
\|\omega(\cdot,y)\|_{L^2(\Omega_0)} \le C \|\omega(\cdot,y)\|_{W^{2,\frac{2n}{n+4}}(\Omega_0)}
\\
\le C
\|(q_m-q)H(\cdot,y)\|_{L^{\frac{2n}{n+4}}(\Omega_0)}
\le
\left\{
\begin{aligned}
  &\quad \quad C &&, \quad & 4 \le n\le 7 ,
  \\
  &|\ln(\operatorname{dist}(y, \partial D_m))| &&, \quad & n=8 ,
  \\
  &{\operatorname{dist}(y, \partial D_m)}^{4- \frac{n}{2}} &&, \quad & n\ge 9.
\end{aligned}
\right.
\end{multline}
Then, using the asymptotic behavior of the Hankel function again and the inequality
\[
\|G\|_{L^2(\Omega_0 \backslash B_r(y))} \le \|\omega\|_{L^2(\Omega_0 \backslash B_r(y))} + \|H\|_{L^2(\Omega_0 \backslash B_r(y))}
\]
we immediately get (\ref{singular_func_est}).


Let $\tilde{\Gamma}(\cdot)$ stand for the Gamma function.
Noting that
\[
H(\cdot,y) \, , \, \Gamma(\cdot, y) \in C^\infty(\Omega_0 -\{y\})
\]
and
\[
\begin{aligned}
& H(x,y) - \Gamma(x,y)
\\
\sim \,\, & -\frac{\ii}{\pi}\tilde{\Gamma}\left(\frac{n-2}{2}\right)\frac{1}{4 \ii \, \pi^{(n-2)/2}}|x-y|^{2-n} - \frac{{\tilde{\Gamma}\left(\frac{n+2}{2}\right)}}{n(2-n)\pi^{n/2}}|x-y|^{2-n}
\\
= \,\, & 0 ,
\end{aligned}
\]
as $|x-y|$ goes to $0$, we conclude that $|\Gamma(\cdot,y) - H(\cdot,y)|$ is uniformly bounded for all $y$ such that $\operatorname{dist}(y,\partial \Omega_0)\ge \frac{r_0}{16}$. Then 
(\ref{est_green2}) follows.
\end{proof}

\medskip\medskip

In both Beretta \& Francini's
proof \cite{Beretta2011} and Alessandrini \& Vessella's proof
\cite{Alessandrini2005}, the blow-up property of a singular function,
\[
   \int_{U_k} \nabla G_1(y,x) \nabla G_2(x,y) \, \dd x ,
\]
where $U_k = \Omega \backslash \cup_{j = 1}^k D_j$ and $G_1$, $G_2$
are functions defined by (\ref{fS}) for potentials $q^{(1)}$,
$q^{(2)}$, respectively, when $y$ approaches the interfaces, is
essential. However, in the case of the Schr\"{o}dinger equation, this
does not happen if $n = 2,3$. Therefore, for $n = 2,3$, we will
introduce a derivative in the point source. For $n=3$, let
\[
   \Gamma(x,y) = -\frac{x_3-y_3}{4 \pi |x-y|^3} ,
\]
which is the solution to the equation
\begin{equation}
   -\Delta \Gamma(x,y) =
           \frac{\partial}{\partial x_3} \delta_y(x) .
\end{equation}

\medskip\medskip

\begin{proposition}\label{G_behavior-3D}
Let $n = 3$ and $q\in L^\infty(\Omega_0)$. For $y\in \Omega_0$, there exists a unique function $G(\cdot, y)$ continuous in $\Omega_0 \backslash \{ y \}$ such that
\begin{equation}
\int_{\Omega_0} \nabla G(\cdot, y) \cdot \nabla \phi +q G(\cdot, y)\phi =  \frac{\partial}{\partial x_n}  \phi(y), \quad \forall \phi\in C_0^\infty(\Omega).
\end{equation}
Furthermore, we have that $G(x,y)$ is symmetric, i.e.,
\begin{equation}\label{symmetry3d}
G(x,y) = G(y,x), \quad x,y\in \Omega_0,
\end{equation}
and the following estimates
\begin{equation}\label{singular_func_est-3d}
\|G(\cdot, y)\|_{L^2(\Omega_0 \backslash B_r(y))} \le Cr^{-\frac{1}{2}}, \quad r\le \frac{1}{2} \mbox{dist} (y, \partial \Omega)
\end{equation}
and
\begin{equation}\label{est_green1-3d}
\|G(\cdot, y) - \Gamma(\cdot,y)\|_{L^2(\Omega_0)} \le C, \quad \mbox{dist}(y, \partial \Omega_0) \ge \frac{r_0}{16}.
\end{equation}
hold, where the constant $C$ depends on the constant in Proposition~\ref{2-energy}.
\end{proposition}

\medskip\medskip

\begin{proof}
Consider $G(x,y) = \Gamma(x,y) +\omega (x,y)$, where $\omega$ solves
\begin{equation}\label{remainder-equ}
\left\{
\begin{array}{rl}
(-\Delta + q)\omega = & q \Gamma, \quad \mbox{ in }\Omega_0, \\
\omega = & - \Gamma, \quad \mbox{ on }\partial \Omega_0.
\end{array}
\right.
\end{equation}
Since $\Gamma(\cdot, y) \in W^{5/4, 4/3}(\partial \Omega_0)$, $q\Gamma \in L^{4/3}(\Omega_0)$ and $-\Gamma(\cdot,y)\in H^{1/2}(\partial\Omega_0)$, by Proposition~$\ref{2-energy}$,   $(\ref{remainder-equ})$ has a unique solution $\omega \in H^{1}( \Omega_0)$ and $\omega = G - \Gamma$ satisfies the estimate
\begin{equation}\label{omegaestimate}
\| \omega(\cdot, y) \|_{H^1(\Omega_0)} \le C \left( \| \Gamma(\cdot, y) \|_{H^{1/2}(\partial \Omega_0)} + \| q(\cdot)\Gamma(\cdot, y) \|_{H^{-1}(\Omega_0)} \right) \le C,
\end{equation}
when $\operatorname{dist}(y,\partial \Omega_0) \ge \frac{r_0}{16}$. Hence,
\begin{multline}\label{res-est}
 \| G(\cdot, y)-\Gamma(\cdot, y) \|_{L^2(\Omega_0\backslash B_r(y))} = \| \omega(\cdot, y) \|_{L^2(\Omega_0\backslash B_r(y))} \\
 \le \| \omega(\cdot, y) \|_{H^1(\Omega_0\backslash B_r(y))}\le \| \omega(\cdot, y) \|_{H^1(\Omega_0)} \le C.
\end{multline}
With the fact that
\begin{equation}
\|\Gamma(\cdot,y)\|_{L^2(\Omega_0 \backslash B_r(y))} \le Cr^{-\frac{1}{2}}, \quad r\le \frac{1}{2} \mbox{dist} (y, \partial \Omega),
\end{equation}
(\ref{res-est}) gives the desired estimate
\begin{equation}
\|G(\cdot,y)\|_{L^2(\Omega_0 \backslash B_r(y))} \le Cr^{-\frac{1}{2}}, \quad r\le \frac{1}{2} \mbox{dist} (y, \partial \Omega).
\end{equation}

Finally, again  by (\ref{omegaestimate}) we immediately get
\begin{equation}
\| G(\cdot, y)-\Gamma(\cdot, y) \|_{L^2(\Omega_0)}  \le C.
\end{equation}
\end{proof}

\medskip\medskip

For $n=2$, let
\[
   \Gamma(x,y) = -\frac{2 \pi (x_2-y_2)}{|x-y|^2 }
\]
which is the solution to the equation
\begin{equation}
   -\Delta \Gamma(x,y) =
       \frac{\partial}{\partial x_2} \delta_y(x) .
\end{equation}

\medskip\medskip

\begin{proposition}\label{G_behavior-2D}
Let $n = 2$ and $q \in L^\infty(\Omega_0)$. For $y\in \Omega_0$, there
exists a unique function $G(\cdot, y)$ continuous in $\Omega_0
\backslash \{ y \}$ such that
\begin{equation}
   \int_{\Omega_0} (\nabla G(\cdot,y) \cdot \nabla \phi +
      q G(\cdot,y) \phi
   = \frac{\partial}{\partial x_n} \phi(y) ,\quad
     \forall \phi \in C_0^\infty(\Omega) .
\end{equation}
Furthermore, we have that $G(x,y)$ is symmetric, that is,
\begin{equation}\label{symmetry2d}
   G(x,y) = G(y,x) ,\quad x,y \in \Omega_0 ,
\end{equation}
and the estimates
\begin{equation}\label{singular_func_est-2d}
   \| G(\cdot,y) \|_{L^2(\Omega_0 \backslash B_r(y))}
        \le C |\ln{r}|^{\frac{1}{2}} ,\quad
   r \le \frac{1}{2}
     \min\left(
     \operatorname{dist}(y,\partial \Omega_0),\frac{1}{2}\right)
\end{equation}
and
\begin{equation}\label{est_green1-2d}
   \| G(\cdot,y) - \Gamma(\cdot,y) \|_{L^2(\Omega_0)} \le C ,\quad
   \operatorname{dist}(y,\partial \Omega_0) \ge \frac{r_0}{16}
\end{equation}
hold, where the constant $C$ depends on the constant in Proposition~\ref{2-energy}.
\end{proposition}

\medskip\medskip

We omit the proof here, because it follows from an adaption of the proof
of Proposition~\ref{G_behavior-3D}. The symmetry of $G$ follows by
standard arguments based on integration by parts (see for example
\cite{Evans2010}).

\medskip\medskip
In the sequel we will derive estimates of unique continuation in $K$ for solutions to our equation.
A key ingredient to obtain these estimates is the Three Spheres Inequality that we will state below and that was proved by \cite[Theorem~3.1]{Alessandrini2001}.
The next two propositions concern Three Sphere Inequalities for our
equation. To prove it, one interprets the equation $(-\Delta +q) u = 0$ for a complex function $q(x)$ as a weakly coupled  system of equations with Laplacian
principal part
\begin{equation}\label{weaklysystem}
-\Delta U+QU=0,
\end{equation}
where $U$ is a vector with components the real and imaginary parts of
$u$, that is, $u^{(1)} = \mathfrak{R} u$, $u^{(2)} =\mathfrak{I} u$,
and $Q$ is a two by two tensor with elements the real and complex part
of the potential $q$, that is, $q^{(1)} = \mathfrak{R} q$ and $
q^{(2)} =\mathfrak{I} q$. We can also write the system in the form
\begin{equation*}
\left\{
\begin{array}{rrl}
   -\Delta u^{(1)} + q^{(1)}u^{(1)} - q^{(2)}u^{(2)} & = & 0 ,
\\[0.2cm]
   -\Delta u^{(2)} + q^{(1)}u^{(2)} + q^{(2)}u^{(1)} & = & 0 .
\end{array}
\right.
\end{equation*}
In \cite[Theorem~3.1]{Alessandrini2001} the authors prove the validity of the Three spheres inequality for elliptic systems with Laplacian principal part. In particular it applies to solutions $U$ of  (\ref{weaklysystem}) and hence also to solutions of $(-\Delta +q) u = 0$.

\medskip\medskip

\begin{proposition}\label{three-sphere-inequality-L2}
Let $u$ be a solution to the equation
\begin{equation*}
   (-\Delta + q) u = 0 \quad \mbox{ in } B_R .
\end{equation*}
Then, for every $\rho_1, \rho_2, \rho_3$, with $0< \rho_1 < \rho_2 <
\rho_3 \le R$,
\begin{equation}
   \| u \|_{L^2(B_{\rho_2})} \le Q_2
        \| u \|_{L^2(B_{\rho_1})}^\alpha
                 \| u \|_{L^2(B_{\rho_3})}^{1- \alpha} ,
\end{equation}
where $\alpha =
\frac{\ln{\frac{\rho_3}{\rho_2}}}{\ln\frac{\rho_3}{\rho_1}} \in
(0,1)$ and $Q_2 \ge 1$ depends on $\|q\|_{L^\infty(B_R)}$ ,
$\frac{\rho_2}{\rho_1}$ and $\frac{\rho_3}{\rho_2}$.
\end{proposition}

\medskip\medskip

\begin{remark}
  In \cite[Theorem~3.1]{Alessandrini2001} the authors prove the validity of the three-spheres inequality for elliptic systems with some limitations on the radii. The derivation of the inequality for arbitrary radii follows by applying the argument of the proof of \cite[Theorem~5.1]{Alessandrini2009} choosing $B_{r_0}(x_0) = B_{r_1}$, $G = B_{r_2}$ and $\Omega = B_{r_3}$.
\end{remark}

\medskip\medskip

Also, we have

\medskip\medskip

\begin{corollary}\label{three-sphere-inequality-L-infty}
Let $u$ be a solution to the equation
\begin{equation*}
   (-\Delta + q) u = 0 \quad \mbox{ in } B_R .
\end{equation*}
Then, for every $\rho_1, \rho_2, \rho_3$, with $0< \rho_1 < \rho_2 <
\rho_3 \le R$,
\begin{equation}
   \| u \|_{L^\infty(B_{\rho_2})} \le Q_\infty
     \| u \|_{L^\infty(B_{\rho_1})}^\beta
            \| u \|_{L^\infty(B_{\rho_3})}^{1- \beta} ,
\end{equation}
where $\beta = \frac{\ln{\frac{2\rho_3}{\rho_2 +
      \rho_3}}}{\ln\frac{\rho_3}{\rho_1}} \in (0,1)$ and $Q_\infty
\ge 1$ depends on $\|q\|_{L^\infty(B_R)}$, $\frac{\rho_2}{\rho_1}$ and
$\frac{\rho_3}{\rho_2}$.
\end{corollary}

\medskip\medskip

\begin{proof}
We use the local boundedness estimate for $u^{(1)}$ and $u^{(2)}$, weak solutions of  elliptic
equations (see for instance \cite[Theorem 8.17]{Gilbarg1983}), to
obtain that there exists a constant $C$, which only depends on $n$ and
$\|q\|_{L^\infty(B_R)}$, such that
\begin{equation}
   \| u \|_{L^{\infty}(B_{\rho_2})} \le
        \frac{C}{(\rho_3 - \rho_2)^{n/2}} \| u \|_{L^2(B_{\rho_3})} .
\end{equation}
Then, by Proposition~\ref{three-sphere-inequality-L2},
\begin{equation}
\begin{aligned}
   \| u \|_{L^\infty(B_{\rho_2})} \le &
     \frac{C}{\left(\frac{\rho_2 + \rho_3}{2} - \rho_2\right)^{n/2}}
           \| u \|_{L^2(B_{\frac{\rho_2 + \rho_3}{2}})}
\\
   \le & \frac{C Q_2}{
    \left(\frac{\rho_2 + \rho_3}{2} - \rho_2\right)^{n/2}}
    \| u \|_{L^2(B_{\rho_1})}^{\alpha}
           \| u \|_{L^2(B_{\rho_3})}^{1- \alpha}
\\
   \le & \frac{C Q_2}{
    \left(\frac{\rho_2 + \rho_3}{2} - \rho_2\right)^{n/2}}
    |B_{\rho_1}|^{\alpha /2} |B_{\rho_3}|^{(1 - \alpha)/2}
    \| u \|_{L^\infty(B_{\rho_1})}^{\alpha}
           \| u \|_{L^\infty(B_{\rho_3})}^{1- \alpha} .
\end{aligned}
\end{equation}
\end{proof}

\medskip\medskip

As a consequence of the Three Spheres Inequality stated in
Corollary~\ref{three-sphere-inequality-L-infty}, we derive the
following quantitative estimate for unique continuation of solutions
to our equation.

\medskip\medskip

\begin{proposition}\label{smallness-propagation}
Let $K$ and $K_0$ be defined as before, and let $v \in H^1(K)$ be a
weak solution to the equation
\begin{equation*}
   (-\Delta + q(x)) v = 0 \quad \mbox{ in } K .
\end{equation*}
Assume that, for given positive numbers $\varepsilon_0$, $E_0$ and real number $\gamma$, $v$
satisfies
\begin{equation}\label{K0-bound}
   \| v \|_{L^\infty(K_0)} \le \varepsilon_0 ,
\end{equation}
and
\begin{equation}\label{infty-bound}
   |v(x)| \le (\varepsilon_0 + E_0) \,
      \operatorname{dist}(x,\Sigma_{M+1})^{\gamma} ,\quad
   x \in K .
\end{equation}
Then the following inequality holds true for every $0 < r < 2 r_1$,
\begin{equation}\label{smallness-estimate}
   |v(\tilde{x})|  \le  C
   \left(\frac{\varepsilon_0}{\varepsilon_0 + E_0}
        \right)^{\tau_r\beta^{N_1}}
     (\varepsilon_0 + E_0) r^{(1 - \tau_r)\gamma} ,
\end{equation}
where $\tilde{x} = P_{M+1} - r \nu(P_{M+1})$ with $\nu$ being the exterior unit normal
vector to $\partial D_M$ 
at $P_{M+1}$,
$\beta = \frac{\ln{(8/7)}}{\ln{4}}$, $\tau_r =
\frac{\ln{\left(\frac{12 r_1 - 2 r}{12 r_1 - 3
      r}\right)}}{\ln{\left(\frac{6 r_1 - r}{2r_1}\right)}} \in
(0,1)$ and the constants $N_1$ and $C$ depend on $r_0, L, A, B$ and $n$.
\end{proposition}

\medskip\medskip

\begin{proof}
We construct a chain of spheres of radius $r_1$
with centers $x_0, x_1, \dots, x_k$ such that the first is
$B_{r_1}(x_0) \subset B_{4 r_1}(x_0) \subset K_0$, all the spheres are
externally tangent, and the last one is centered at $x_k = P_{M + 1} -
3 r_1 \nu (P_{M + 1})$. We choose this chain so that the spheres of
radius $4r_1$ concentric with those of the chain, except the last one,
are contained in $K$ and have a distance greater than $r_1$ away from
$\Sigma_{M + 1}$. Such a chain has a finite number of spheres that is
smaller than $N_1 = \frac{A}{|B_{r_1}|} + 1$.

By Corollary~\ref{three-sphere-inequality-L-infty} and
(\ref{infty-bound}), we have
\begin{equation*}
\begin{array}{rl}
   \| v \|_{L^\infty(B_{r_1}(x_1))}
   \le & \|v\|_{L^\infty(B_{3r_1}(x_0))}
\\[0.3cm]
   \le & Q_\infty \| v \|_{L^\infty(B_{r_1}(x_0))}^\beta
           \| v \|_{L^\infty(B_{4r_1}(x_0))}^{1 - \beta}
\\
   \le & \displaystyle{
   C \left(\frac{\varepsilon_0}{\varepsilon_0 + E_0}\right)^{\beta}
              (\varepsilon_0 + E_0)} ,
\end{array}
\end{equation*}
where $C$ depends on $Q_\infty$ and $r_1$. By iterated application of
Corollary~\ref{three-sphere-inequality-L-infty} to $v$ with radii
$r_1$, $3 r_1$ and $4 r_1$ over the chain of spheres, we have, by
(\ref{K0-bound}),
\begin{equation*}
\begin{array}{rl}
   \| v \|_{L^\infty(B_{r_1}(x_k))}
   \le & Q_\infty \| v \|_{L^\infty(B_{r_1}(x_{k-1}))}^\beta
         \| v \|_{L^\infty(B_{4r_1}(x_{k-1}))}^{1 - \beta}
\\
   \le & \displaystyle{
   C \left(\frac{\varepsilon_0}{\varepsilon_0
       + E_0}\right)^{\beta^{N_1}} (\varepsilon_0 + E_0)} ,
\end{array}
\end{equation*}
where $C$ depends on $Q_\infty$ and $r_1$. Now, we let $\tilde{x} =
P_{M+1} - r \nu(P_{M+1})$ where $r < 2 r_1$. Using
Corollary~$\ref{three-sphere-inequality-L-infty}$ again for spheres
centered at $x_k$ of radii $r_1$, $3r_1 - r$ and $3r_1 - \frac{r}{2}$,
we obtain that
\begin{equation*}
\begin{array}{rl}
   \| v \|_{L^\infty(B_{3r_1 - r}(x_k))}
   \le & Q_\infty \| v \|_{L^\infty(B_{r_1}(x_k))}^{\tau_r}
         \| v \|_{L^\infty(B_{3r_1 - \frac{r}{2}}(x_k))}^{1 - \tau_r}
\\
   \le & \displaystyle{
   C \left(\frac{\varepsilon_0}{\varepsilon_0
       + E_0}\right)^{\tau_r\beta^{N_1}} (\varepsilon_0 + E_0)
              r^{(1 - \tau_r)\gamma}} ,
\end{array}
\end{equation*}
which completes the proof.
\end{proof}

\medskip\medskip

\begin{remark}\label{smallness-propagation-2D}
Let us observe that, in order to apply Proposition~\ref{smallness-propagation} to the singular function defined in Section~\ref{sec:4} when $n=2,4$, we need to replace the condition (\ref{infty-bound}) by
\begin{equation}
   |v(x)| \le (\varepsilon_0 + E_0) \,
      |\ln (\operatorname{dist}(x,\Sigma_{M+1}))|^{\frac{1}{2}} ,\quad
   x \in K.
\end{equation}
By using the same proof technique, we can obtain the same result with (\ref{smallness-estimate}) replaced by
\begin{equation}
   |v(\tilde{x})|  \le  C
   \left(\frac{\varepsilon_0}{\varepsilon_0 + E_0}
        \right)^{\tau_r\beta^{N_1}}
     (\varepsilon_0 + E_0) |\ln r|^{\frac{1 - \tau_r}{2}} .
\end{equation}
\end{remark}

\section{Proof of the main result}

\label{sec:4}

Assume that $D_{M}$ is the subdomain of the partition of $\Omega$  where the maximum of $\|q^{(1)} - q^{(2)}\|$ is realized and let us denote
\begin{equation}
   E = \| q^{(1)} - q^{(2)}\|_{L^\infty(D_{M})}
              = \| q^{(1)} - q^{(2)}\|_{L^\infty(\Omega)}.
\end{equation}
We consider the chain of domains, $D_0,D_1,\dots,D_{M}$, as before; $S,K$ and $K_0$ are defined as in the previous section. We set
\[
\displaystyle{U_0 = \Omega, \, U_k = \Omega \backslash \cup_{j = 1}^k D_j, \quad k = 1, \dots, M \mbox{ and }W_k = \cup_{j = 0}^k D_j.}
\]
Let $y \in K$. For dimension $n \ge 4$, let $G_1(x,y)$ and $G_2(x,y)$ be the Green's function related to $q^{(1)}$ and $q^{(2)}$, respectively, the existence and behavior of
which was shown in Proposition~\ref{G_behavior}. For dimension $n = 2,3$, let $G_1(x,y)$ and $G_2(x,y)$ be a first order derivative of the Green's function, the existence and behavior of which was shown in Propositions~\ref{G_behavior-2D} and \ref{G_behavior-3D}, respectively. We define
\begin{equation}\label{singular_func_def}
   S_k(y,z) = \int_{U_k} (q^{(1)} - q^{(2)})(x)
                           G_1(x,y) G_2(x,z) \, \dd x .
\end{equation}
By Proposition~\ref{G_behavior}, \ref{G_behavior-3D} and \ref{G_behavior-2D}, there exist a
constant $C$ such that
\begin{equation}
\begin{aligned}
   |S_k(y,z)| \le & C E\
   |\ln(\operatorname{dist}(y, U_k)) \,
      \ln(\operatorname{dist}(z, U_k))|^{\frac{1}{2}} ,\quad
   & y,z \in K \cap W_k , \quad & n=2,4;
   \\
   |S_k(y,z)| \le & C E\
   (\operatorname{dist}(y, U_k) \,
      \operatorname{dist}(z, U_k))^{-\frac{1}{2}} ,\quad
   & y,z \in K \cap W_k , \quad & n=3;
   \\
   |S_k(y,z)| \le & C E\
   (\operatorname{dist}(y, U_k) \,
      \operatorname{dist}(z, U_k))^{2-\frac{n}{2}} ,\quad
   & y,z \in K \cap W_k , \quad & n\ge 5.
   \end{aligned}
\end{equation}
We focus on $n = 3$ first; we will discuss the adaptation of the proof for the case $n = 2,4$ and $n \ge 5$ at the end of the proof.

\medskip\medskip

\begin{lemma}\label{prop-sol-singular}
For every $y,z \in K \cap W_k$, we have $S_k(\cdot, z)$, $S_k(y,
\cdot)\in H^1(K\cap W_k)$ and
\begin{equation}
   (-\Delta + q^{(1)}) S_k(\cdot,z) = 0 ,\quad
   (-\Delta + q^{(2)}) S_k(y,\cdot) = 0  \quad
   \mbox{ in } K\cap W_k .
\end{equation}
\end{lemma}

\medskip\medskip

\noindent

The proof of this Lemma follows from the symmetry of $G_i\ (i = 1,2)$
and changing the order of integration and differentiation.\\

\medskip\medskip

\begin{lemma}\label{Prop-Sk-smallness}
If for some $\varepsilon_0 > 0$ and $k \in \{1, \dots, M-1\}$ we have
that
\begin{equation}
   |S_k(y,z)| \le \varepsilon_0 ,\quad \forall y,z \in K_0 ,
\end{equation}
then
\begin{equation}\label{smallness-propagation-Sk}
\begin{aligned}
   |S_k(y_r,y_r)| \le & C
   \left(\frac{\varepsilon_0}{\varepsilon_0
              + E}\right)^{\tau_r^2 \beta^{2N_1}} \,
   (\varepsilon_0 + E) \, |\ln r| , \quad && n = 2 \mbox{ or } 4,
   \\
   |S_k(y_r,y_r)| \le & C
   \left(\frac{\varepsilon_0}{\varepsilon_0
              + E}\right)^{\tau_r^2 \beta^{2N_1}} \,
   (\varepsilon_0 + E) \, r^{-1} , \quad && n = 3,
   \\
   |S_k(y_r,y_r)| \le & C
   \left(\frac{\varepsilon_0}{\varepsilon_0
              + E}\right)^{\tau_r^2 \beta^{2N_1}} \,
   (\varepsilon_0 + E) \, r^{4-n} , \quad && n \ge 5,
   \end{aligned}
\end{equation}
where $y_r = P_{k+1} - r \nu(P_{k+1})$, $r$ is small,
$\nu(P_{k+1})$ is the exterior unit normal vector to $\partial D_k$ at
$P_{k+1}$ and the positive constant $C$ depends on $r_0, L, A, B$ and $n$.
\end{lemma}

\medskip\medskip

\begin{proof}
Let the dimension $n=3$. We fix $z\in K_0$ first and consider $v(y) = S_k(y,z)$. By
Lemma~\ref{prop-sol-singular}, $v$ solves the equation $(-\Delta +
q^{(1)})v = 0$ in $K\cap W_k.$ Moreover, by (\ref{singular_func_est-3d}),
we have
\begin{equation}
   |v(y)| \le C \, E\
   \operatorname{dist}(y,\Sigma_{k + 1})^{-\frac{1}{2}} ,\quad
   y \in K \cap W_k .
\end{equation}
Then, by Proposition~\ref{smallness-propagation} with $\gamma = -\frac{1}{2}$, we have, for $0 < r < 2 r_1$,
\begin{equation}
   |S_k(y_r,z)| \le C
   \left(\frac{\varepsilon_0}{\varepsilon_0
                + E}\right)^{\tau_r\beta^{N_1}}
   (\varepsilon_0 + E) \, r^{-\frac{1}{2}} .
\end{equation}
Next, we consider
\begin{equation}
   \tilde{v}(z) = S_k(y_r,z) ,\quad z \in K\cap W_k ,
\end{equation}
which solves the equation $(-\Delta + q^{(2)}) \tilde{v} = 0$ in $K
\cap W_k$, and, by (\ref{singular_func_est-3d}), satisfies
\begin{equation}
   |\tilde{v}(z)| \le C \, E\
   \left( r \operatorname{dist}(z,\Sigma_{k + 1})
          \right)^{-\frac{1}{2}} ,\quad z \in K\cap W_k .
\end{equation}

By Proposition~\ref{smallness-propagation}, again, we then obtain
estimate (\ref{smallness-propagation-Sk}) for $n=3$.

The proof for other dimensions follows from the same proof with a few modifications. For $n=2,4$, a modified version of Proposition~\ref{smallness-propagation}, as stated in Remark~\ref{smallness-propagation-2D}, needs to be applied. For $n \ge 5$, one can apply Proposition~\ref{smallness-propagation} with $\gamma = 2 -\frac{n}{2}$.
\end{proof}

\medskip\medskip

\noindent

\textit{Proof of Theorem~\ref{thm-stability}}. Let
\[
\varepsilon = \|\Lambda_1^{(\Sigma)} - \Lambda_2^{(\Sigma)}\|_{\mathcal{L}(H^{1/2}, H^{-1/2})}
\]
 and
\[
 \delta_k = \|q^{(1)} - q^{(2)}\|_{L^\infty(W_k)}, \quad k = 0,1,\dots, M.
 \]
From
the Alessandrini identity (see for instance, Chapter $5$ of
\cite{Isakov2006})
\begin{equation}
   \int_{\Omega} (q^{(1)} - q^{(2)})(x)
      G_1(x,y) G_2(x,z) \, \dd x  =
   \langle (\Lambda_1 - \Lambda_2) G_1(\cdot,y),
      \overline{G}_2(\cdot,z) \rangle ,\quad \forall y,z \in K_0
\end{equation}
and Proposition~\ref{G_behavior}, we find that
\begin{equation}\label{Sk-K0}
   |S_{k-1}(y,z)| \le C \, (\varepsilon + \delta_{k-1}) .
\end{equation}
Let $P_k \in \Sigma_k$ and $y_r = z_r = P_k - r \nu(P_k)$, where $\nu(P_k)$ is the exterior unit normal vector to
$\partial D_{k-1}$ 
and $r$ is small. We write
\begin{equation}\label{dec-S}
   S_{k-1}(y_r,y_r) = I_1 + I_2
\end{equation}
with
\begin{equation}\label{def-I1}
   I_1 = \int_{B_{\rho_0}(P_k)\cap D_k}
           (q^{(1)} - q^{(2)})(x) G_1(x,y_r) G_2(x,y_r) \, \dd x
\end{equation}
and
\begin{equation}\label{def-I2}
   I_2 = \int_{U_{k-1} \backslash (B_{\rho_0}(P_k) \cap D_k)}
           (q^{(1)} - q^{(2)})(x) G_1(x,y_r) G_2(x,y_r) \, \dd x ,
\end{equation}
where $\rho_0 = \frac{r_0}{6}.$

For $n = 3$, by Proposition~\ref{G_behavior-3D}, we have
\begin{equation}\label{est-I2}
   |I_2| \le C \, E .
\end{equation}
We estimate $I_1$ as follows:
\begin{equation}\label{est-I1}
  \begin{aligned}
   |I_1|  = & |q^{(1)}_k - q^{(2)}_k| \,
     \left| \int_{B_{\rho_0}(P_k) \cap D_k} \!\!
              G_1(x,y_r) G_2(x,y_r) \dd x \right|
\\
   \ge & |q^{(1)}_k - q^{(2)}_k| \left\{\mstrut{0.5cm}\right.
       \left| \int_{B_{\rho_0}(P_k) \cap D_k} \!\!
       \Gamma(x,y_r) \Gamma(x,y_r) \dd x \right|
\\
   &- \left| \int_{B_{\rho_0}(P_k) \cap D_k} \!\!
     (G_1(x,y_r) - \Gamma(x,y_r)) \Gamma(x,y_r) \dd x \right|
\hspace*{3.0cm}
\\
\hspace*{1.0cm}
  & - \left| \int_{B_{\rho_0}(P_k) \cap D_k} \!\!
     (G_2(x,y_r) - \Gamma(x,y_r))\Gamma(x,y_r) \dd x \right|
\\
   & - \left| \int_{B_{\rho_0}(P_k) \cap D_k} \!\!
     (G_1(x,y_r) - \Gamma(x,y_r))
     (G_2(x,y_r) - \Gamma(x,y_r)) \dd x \right|
     \left.\mstrut{0.5cm}\right\} .
\end{aligned}
\end{equation}
By Propositions~\ref{G_behavior-3D} and the fact that
\begin{multline*}
   \left| \int_{B_{\rho_0}(P_k)\cap D_k} \!\!
          (G_i(x,y_r) - \Gamma(x,y_r))
           \Gamma(x,y_r) \dd x \right|
\\
   \le \frac{1}{2} \int_{B_{\rho_0}(P_k)\cap D_k} \!\!
   \left( 2 |G_i(x,y_r) - \Gamma(x,y_r)|^2
          + \frac{1}{2} |\Gamma(x,y_r)|^2 \right) \dd x ,\quad  i = 1,2
\end{multline*}
and
\begin{multline*}
    \left| \int_{B_{\rho_0}(P_k) \cap D_k} \!\!
     (G_1(x,y_r) - \Gamma(x,y_r))
     (G_2(x,y_r) - \Gamma(x,y_r)) \dd x \right|
     \left.\mstrut{0.5cm}\right\}
\\
   \le \frac{1}{2} \int_{B_{\rho_0}(P_k)\cap D_k} \!\!
   \left( |G_1(x,y_r) - \Gamma(x,y_r)|^2
          + |G_2(x,y_r) - \Gamma(x,y_r)|^2 \right) \dd x ,
\end{multline*}
we obtain that
\[
  |I_1| \ge  |q^{(1)}_k - q^{(2)}_k| \left(\mstrut{0.5cm}\right.
       \frac{1}{2} \int_{B_{\rho_0}(P_k) \cap D_k} \!\!
       |\Gamma(x,y_r)|^2 \dd x - C
     \left.\mstrut{0.5cm}\right) .
\]
Using the explicit form of
$\Gamma(x,y)$, we find that
\begin{eqnarray}
   |I_1| & \ge & |q^{(1)}_k - q^{(2)}_k| (C r^{-1} - C)
\nonumber\\
         & \ge & C |q^{(1)}_k - q^{(2)}_k |r^{-1} - C E .
\label{eq:4-17}
\end{eqnarray}
Now, by Lemma~\ref{Prop-Sk-smallness} and (\ref{Sk-K0}), we have
\begin{equation*}
   |S_{k-1}(y_r, y_r)| \le C \left(
   \frac{\varepsilon + \delta_{k-1}}{\varepsilon + \delta_{k-1}
         + E} \right)^{\tau_r^2 \beta^{2N_1}}
             (\varepsilon + \delta_{k-1} + E) r^{-1} .
\end{equation*}
Hence, using (\ref{dec-S}), (\ref{est-I2}) and (\ref{eq:4-17}), we have
\begin{equation*}
   |q^{(1)}_k - q^{(2)}_k| \, r^{-1}
   \le C \left(\mstrut{0.5cm}\right.
   E + \left(\frac{\varepsilon + \delta_{k-1}}{\varepsilon
       + \delta_{k-1} + E}\right)^{\tau_r^2\beta^{2N_1}}
       (\varepsilon + \delta_{k-1} + E) r^{-1}
       \left.\mstrut{0.5cm}\right) ,
\end{equation*}
so that
\begin{equation}
   |q^{(1)}_k - q^{(2)}_k|
   \le C (\varepsilon + \delta_{k-1} + E)
   \left(\mstrut{0.5cm}\right.
   \left(\frac{\varepsilon + \delta_{k-1}}{\varepsilon
         + \delta_{k-1} + E}\right)^{\tau_r^2\beta^{2N_1}}
   + r \left.\mstrut{0.5cm}\right) .
\end{equation}

Noting that
\[
\displaystyle{\tau_r =
\frac{\ln{\left(\frac{12 r_1 - 2 r}{12 r_1 - 3
      r}\right)}}{\ln{\left(\frac{6 r_1 - r}{2r_1}\right)}}}, \quad \forall  r \in (0, 2 r_1)
\]
implies
\[
\displaystyle{ \frac{\tau_r}{r} \ge \frac{1}{12r_1 \ln{3}} }, \quad \forall  r \in (0, 2 r_1)
\]
we get
\begin{equation}
   |q^{(1)}_k - q^{(2)}_k|
   \le C (\varepsilon + \delta_{k-1} + E)
   \left(\mstrut{0.5cm}\right.
   \left(\frac{\varepsilon + \delta_{k-1}}{\varepsilon
         + \delta_{k-1} + E}\right)^{\beta^{2N_1} (12 r_1 \ln{3})^{-2} r^2}
   + r \left.\mstrut{0.5cm}\right) .
\end{equation}

By taking $r = \left|\ln{\left(\frac{\varepsilon +
    \delta_{k-1}}{\varepsilon + \delta_{k-1} +
    E}\right)}\right|^{-1/4}$ and noting that
    \[
    \displaystyle{
    \left(e^{-r^{-4}}\right)^{\beta^{2N_1} (12 r_1 \ln{3})^{-2} r^2} \le C r, \quad \forall r>0
    }
    \]
    for some constant $C$, we obtain that
\begin{equation}\label{induc-equ}
   |q^{(1)}_k - q^{(2)}_k|
   \le C \, (\varepsilon + \delta_{k-1} + E)
   \left|\ln{\left(\frac{\varepsilon + \delta_{k-1}}{
        \varepsilon + \delta_{k-1} + E}\right)}\right|^{-\frac{1}{4}} .
\end{equation}

We let
\[
   \omega(t) = \left\{
   \begin{array}{rl}
   |\ln{t}|^{-\frac{1}{4}} , & \quad 0 < t < e^{-3} ,
   \\[0.1cm]
   3^{-\frac{1}{4}} , & \quad t \ge e^{-3} .
   \end{array}\right.
\]
Noting that the function $t \mapsto t \omega_n(1/t)$ is increasing, we have
\[
\frac{\varepsilon + \delta_{k-1} + E}{\varepsilon + \delta_{k-1}} \, \omega \! \left(\frac{\varepsilon + \delta_{k-1}}{
        \varepsilon + \delta_{k-1} + E}\right) \ge \omega(1),
\]
hence
\[
\delta_{k-1} \le \varepsilon + \delta_{k-1} \le (\omega(1))^{-1}(\varepsilon + \delta_{k-1} + E) \,\omega \! \left(\frac{\varepsilon + \delta_{k-1}}{
        \varepsilon + \delta_{k-1} + E}\right),
\]
which with (\ref{induc-equ}) gives that
\begin{equation}\label{induc-equ2}
   \delta_k \le C \, (\varepsilon + \delta_{k-1} + E) \,
   \omega \! \left(\frac{\varepsilon
   + \delta_{k-1}}{\varepsilon + \delta_{k-1} + E}\right) .
\end{equation}
The above choice of $r$ is possible only if
\[
   \left|\ln{\left(\frac{\varepsilon + \delta_{k-1}}{
       \varepsilon + \delta_{k-1} + E}\right)}\right|^{-1/4} < 2 r_1.
\]
However, if
\[
   \left|\ln{\left(\frac{\varepsilon + \delta_{k-1}}{
       \varepsilon + \delta_{k-1} + E}\right)}\right|^{-1/4} \ge 2 r_1,
\]
that is,
\[
\frac{\varepsilon + \delta_{k-1}}{\varepsilon + \delta_{k-1} + E} \ge e^{-(2 r_1)^{-4}},
\]
the fact that
\[
\displaystyle{\sup_{\substack{r \in (0, \, 2r_1) \\ t \in ( e^{-(2 r_1)^{-4}}, \, 1)}} t^{\beta^{2N_1} (12 r_1 \ln{3})^{-2} r^2}|\ln{t}|^{\frac{1}{4}}}
\]
is finite shows that (\ref{induc-equ}) still holds true, then (\ref{induc-equ2}) follows.

We iterate (\ref{induc-equ2}), starting from $\delta_0 = 0$, and find
\begin{equation}
   \delta_k + \varepsilon \le (C + 3^{1/4})^k
   (E + \varepsilon) \,
   \omega_k \! \left(\frac{\varepsilon}{\varepsilon + E}\right) ,
\end{equation}
where $\omega_k$ is the composition of $\omega$ $k$ times with
itself. We recall that $E = \delta_M$, whence,
\begin{equation}\label{fi-equa}
   E + \varepsilon \le (C + 3^{1/4})^M (E + \varepsilon) \,
   \omega_M \! \left(\frac{\varepsilon}{\varepsilon + E}\right) ,
\end{equation}
so that
\begin{equation}
   E \le \frac{1 - \omega^{-1}_M((C + 3^{1/4})^{-M})}{
         \omega^{-1}_M((C + 3^{1/4})^{-M})} \, \varepsilon ,
\end{equation}
which completes the proof for dimension $n = 3$.

The proof for $n = 2$ and $n = 4$ follows from a careful inspection and adaptation
of the above proof for $n = 3$. By Proposition~\ref{G_behavior-2D} and \ref{G_behavior-3D}, and the explicit form of
\[
\begin{aligned}
  \Gamma(x,y) = & -\frac{2 \pi (x_2-y_2)}{|x-y|^2 }, \quad & n=2,
  \\
  \Gamma(x,y) = & -\frac{1}{4 \pi^2 |x-y|^2 }, \quad & n=4,
\end{aligned}
\]
we obtain that
\begin{equation}
   |q^{(1)}_k - q^{(2)}_k| \le C (\varepsilon + \delta_{k-1} + E)
   \left(\mstrut{0.5cm}\right.
   \left(\frac{\varepsilon + \delta_{k-1}}{
   \varepsilon + \delta_{k-1} + E}\right)^{\tau_r^2\beta^{2N_1}}
           + |\ln r|^{-1} \left.\mstrut{0.5cm}\right) .
\end{equation}
Then, by taking $r = \frac{\varepsilon + \delta_{k-1}}{\varepsilon +
  \delta_{k-1} + E}$ and adapting the function $\omega(t)$ according
to
\[
   \tilde{\omega}(t) = \left\{
   \begin{array}{rl}
   |\ln{t}|^{-1} , & \quad 0 < t < e^{-2} ,
   \\[0.1cm]
      \frac{1}{2} , & \quad t \ge e^{-2} ,
   \end{array}\right.
\]
we end up with
\begin{equation}
   E \le \frac{1 - \tilde{\omega}^{-1}_M((C + 2)^{-M})}{
         \tilde{\omega}^{-1}_M((C + 2)^{-M})} \, \varepsilon ,
\end{equation}
which completes the proof for $n=2$ and $n=4$.

Let us sketch the required modifications of the proof for higher dimensional cases($n \ge 5$) below. First, one can use the same decomposition of the singular function as in (\ref{dec-S}), (\ref{def-I1}) and (\ref{def-I2}), and by Proposition~\ref{G_behavior}, the same upper bound estimate of $I_2$ as in (\ref{est-I2}) is obtained. Then, because the order of the singularity of $\Gamma(\cdot,y)$ increases as the dimension $n$ increases, $G(\cdot,y_r) - \Gamma(\cdot,y_r)$ may not be uniformly bounded in $L^2(B_{\rho_0}(P_k)\cap D_k)$ with respect to $r$. A feasible modification here is to compare the orders of singularity of $G(\cdot,y_r) - \Gamma(\cdot,y_r)$ and $\Gamma(\cdot,y_r)$. More precisely, we can estimate
\[
   \left| \int_{B_{\rho_0}(P_k)\cap D_k} \!\!
          (G_i(x,y_r) - \Gamma(x,y_r))
           \Gamma(x,y_r) \dd x \right|,  \quad i=1,2,
\]
using H\"older inequality
, as
\[
\begin{aligned}
   \,\, & \left| \int_{B_{\rho_0}(P_k)\cap D_k} \!\!
          (G_i(x,y_r) - \Gamma(x,y_r))
           \Gamma(x,y_r) \dd x \right|
\\
   \le \,\, & \int_{B_{\rho_0}(P_k)\cap D_k} \!\!
          \left| (G_i(x,y_r) - \Gamma(x,y_r))
           \Gamma(x,y_r)\right| \dd x
\\
\le \,\,& \|G_i(\cdot,y_r) - \Gamma(\cdot,y_r)\|_{L^{2}(\Omega_0)}\, \| \Gamma(\cdot,y_r)\|_{L^{2}(B_{\rho_0}(P_k)\cap D_k)}.
\end{aligned}
\]
Substituting the above inequality into (\ref{est-I1}) and noting the positiveness of $\Gamma(\cdot, y_r)$, we obtain the estimate of the lower bound of $I_1$ as
\begin{equation}\label{est-I1-n}
\begin{aligned}
  |I_1|  \ge   \,\, & |q^{(1)}_k - q^{(2)}_k| ( \|\Gamma(\cdot, y_r)\|_{L^{2}(B_{\rho_0}(P_k)\cap D_k)}^2
  \\
  - \,\,&  \|G_1(\cdot,y_r) - \Gamma(\cdot,y_r)\|_{L^{2}(\Omega_0)}  \|\Gamma(\cdot, y_r)\|_{L^{2}(B_{\rho_0}(P_k)\cap D_k)}
  \\
  - \,\,&  \|G_2(\cdot,y_r) - \Gamma(\cdot,y_r)\|_{L^{2}(\Omega_0)} \|\Gamma(\cdot, y_r)\|_{L^{2}(B_{\rho_0}(P_k)\cap D_k)}
  \\
  - \,\, &  \|G_1(\cdot,y_r) - \Gamma(\cdot,y_r)\|_{L^{2}(\Omega_0)} \, \|G_2(\cdot,y_r) - \Gamma(\cdot,y_r)\|_{L^{2}(\Omega_0)}
  - C).
\end{aligned}
\end{equation}
By the explicit form of
$\Gamma(x,y)$ and Proposition~\ref{G_behavior}, especially (\ref{est_green2}), we observe that $G_i(\cdot,y_r) - \Gamma(\cdot,y_r)$ has the lower order of singularity than $\Gamma(\cdot,y_r)$. 
Hence, by Young's inequality as in the previous proof for $n=3$, we conclude that
\[
\begin{aligned}
  |I_1|  \ge \,\,& C |q^{(1)}_k - q^{(2)}_k |r^{4-n} - C E
\end{aligned}
\]
for $r$ small. Then, following the same argument with the same value of $r$ and noting that
\[
    \displaystyle{
    \left(e^{-r^{-4}}\right)^{\beta^{2N_1} (12 r_1 \ln{3})^{-2} r^2} \le C r^{n-4}, \quad \forall r>0
    }
    \]
    holds true for any $n\ge 5$, we obtain that
\begin{equation}
   |q^{(1)}_k - q^{(2)}_k|
   \le C \, (\varepsilon + \delta_{k-1} + E)
   \left|\ln{\left(\frac{\varepsilon + \delta_{k-1}}{
        \varepsilon + \delta_{k-1} + E}\right)}\right|^{\frac{4-n}{4}} .
\end{equation}
The last step of the modifications is to adapt the function $\omega(t)$ according
to
\[
   \tilde{\omega}(t) = \left\{
   \begin{array}{rl}
   |\ln{t}|^{\frac{4-n}{4}} , & \quad 0 < t < e^{-n} ,
   \\[0.1cm]
      n^{\frac{4-n}{4}} , & \quad t \ge e^{-n} .
   \end{array}\right.
\]
Then we end up with
\begin{equation}
   E \le \frac{1 - \tilde{\omega}^{-1}_M((C + n^{\frac{n-4}{4}})^{-M})}{
         \tilde{\omega}^{-1}_M((C + n^{\frac{n-4}{4}})^{-M})} \, \varepsilon ,
\end{equation}
which completes the proof for $n \ge 5$.
$\square$

\medskip\medskip

\section{Exponential behavior of the Lipschitz stability constant}

\label{sec:5}

In this section, we give a model example to show that the Lipschitz
stability constant $C = C(n, r_0, L, A, N)$ in
Theorem~\ref{thm-stability} behaves exponentially with respect to the
number $N$ of the subdomains. The construction is an analogue of the
construction in \cite{Rondi2006}, pertaining to the inverse
conductivity problem.

Let $\Omega$ be the unit ball $B_1(0)\subset \mathbb{R}^n$ and $D =
[-1/2, 1/2]^n$ be the cube of side $1$ centered at the origin. We
define the class of admissible potentials by
\begin{equation}\label{admissible-set}
   \mathcal{A} = \{ q \in L^\infty(\Omega) \mid 1/2 \le q \le 3/2
   \mbox{ in } \Omega \mbox{ and } q = 1 \mbox{ in }
               \Omega \backslash D \}
\end{equation}
and denote the operator from potential $q$ to $\Lambda_q$ by $F$,
which maps $\mathcal{A}$ into \\ $ \mathcal{L}(H^{1/2}(\partial
\Omega), H^{-1/2}(\partial \Omega))$. We fix a positive integer $N$
and let $N_1$ be the smallest integer such that $N \le N_1^n$. We
divide each side of the cube $D$ into $N_1$ equal parts of length $h =
1/ N_1$ and let $S_{N_1}$ be the set of all open cubes of the type
\[
   D' = (-1/2 + (j_1' - 1)h, -1/2 + j_1' h ) \times \cdots \times
        (-1/2 + (j_n' - 1)h, -1/2 + j_n' h ) ,
\]
where $j_1', \dots,j_n'$ are integers belonging to $\{1, \dots ,
N_1\}$. We order such cubes as follows. For any two different cubes
$D'$ and $D''$ belonging to $S_{N_1}$, we say that $D'\prec D''$ if
and only if there exists an $i_0 \in \{1, \dots, n\}$ such that $j_i'
= j_i''$ for any $i<i_0$ and $j_{i_0}' < j_{i_0}''$. We define
\[
   \mathcal{A}_N = \{ q \in L^\infty(\Omega) \mid
   q(x) = \sum_{j = 1}^N q_j \chi_{D_j}(x) + \chi_{D_0}(x) ,
           \quad q_j \in [1/2, 3/2] \} .
\]
Our aim is to estimate from below the Lipschitz constant $C(N)$ in
terms of $N$. A simple computation shows polynomial behavior of the
lower bound estimate of $C(N)$. To obtain
the exponential estimate, we then need to employ a topological
argument.

Consider a subset $\tilde{\mathcal{A}}_N \subset \mathcal{A}_N$
defined by
\[
   \tilde{\mathcal{A}}_N = \{ q \in L^\infty(\Omega) \mid
   q(x) = \sum_{j = 1}^N q_j \chi_{D_j}(x) + \chi_{D_0}(x) ,
   \quad q_j \in \left\{ \frac{1}{2}, 1 , \frac{3}{2} \} \right\} .
\]
It is easy to check that $\tilde{\mathcal{A}}_N$ is a $1/2$-net of
$\mathcal{A}_N$ with $3^N$ elements and, for any two different $q_1,
q_2 \in \tilde{\mathcal{A}}_N$, we have $\|q_1
-q_2\|_{L^\infty(\Omega)} = 1/2$. Based on Mandache's result
\cite[Lemma~3]{Mandache2001}, there exist a constant $K$, which only
depends on dimension $n$, such that for every $\varepsilon \in (0,
e^{-1})$, there is an $\varepsilon$-net $Y$ for $F(\mathcal{A})$ with
at most $e^{K(-\ln{\varepsilon})^{2n-1}}$ elements. For $\varepsilon \in (0,e^{-1})$ and $N\in\mathbb{N}$ let
\[
Q(\varepsilon, N)= e^{K(-\ln{\varepsilon)}^{2n-1}}.
 \]
 Note that
 \[
 3^N>e^{K(-\ln{\varepsilon})^{2n-1}}
 \]
 if
 \[
 \varepsilon>e^{-K_1N^{1/(2n-1)}}=\varepsilon_0(N)
 \]
 where $K_1 = (K^{-1} \ln 3)^{1/(2n - 1)}$.
There exists $N_0$ such that for $N\geq N_0$ we have that $\varepsilon<e^{-1}$. Thus, for $N\geq N_0$, if we take $\varepsilon=\varepsilon_0$ we have $3^N>Q(\varepsilon,N)$. Then, there exist
two
different $q_1, q_2 \in \tilde{\mathcal{A}}_N$ such that  $\|q_1
-q_2\|_{L^\infty(\Omega)}=1/2$ with their images under
$F$ in the same ball of radius $\varepsilon$ centered at a point of
$Y$, that is,
\[
\frac{1}{2}=\|q_1-q_2\|_{L^\infty(\Omega)}\leq C_N \|\Lambda_{q_1} - \Lambda_{q_2}\|_{
   \mathcal{L}(H^{1/2}(\partial \Omega), H^{-1/2}(\partial \Omega))}\leq 2 C_N\varepsilon_0(N)
\]
from which we get
\[
   C(N) \ge \frac{1}{4}e^{K_1 N^{1/(2n-1)}}.
\]

\section*{Acknowledgements}

This paper was initialized at a Special semester on Inverse Problems
and Applications at MSRI, Berkeley, in the Fall of
2010.  The work of E. Beretta was partially supported by MIUR grant PRIN 20089PWTPS003. The research of M. de Hoop and L. Qiu was supported in part by National Science Foundation grant CMG DMS-1025318, and in part by the members of the Geo-Mathematical Imaging Group at Purdue University.

\def\cprime{$'$} \newcommand{\SortNoop}[1]{}

\bibliographystyle{siam}

\bibliography{conv}

\end{document}